\documentclass{amsart}

\usepackage{amssymb}
\usepackage[centertags]{amsmath}
\usepackage{amsfonts}
\usepackage{epsfig}
\usepackage{amsthm}
\usepackage[ left=40mm, right=40mm, top=30mm, bottom=30mm, includefoot]{geometry}

\newtheorem{conjecture}{Conjecture}[section]
\newtheorem{theorem}{Theorem}[section]
\newtheorem{proposition}[theorem]{Proposition}

\newtheorem{corollary}[theorem]{Corollary}
\theoremstyle{definition}
\newtheorem{definition}[theorem]{Definition}
\theoremstyle{remark}
\newtheorem{remark}[theorem]{Remark}

\newcommand{\Z}{{\mathbb{Z}}}
\newcommand{\s}{{\mathcal{S}}}
\newcommand{\al}{\alpha}
\newcommand{\be}{\beta}
\newcommand{\calL}{{\mathcal{L}}}
\newcommand{\ga}{\gamma}
\newcommand{\BR}{Bollob\'{a}s-Riordan }
\newcommand{\vb}{\textbf{b}}

\begin{document}

\title[Knot invariants and the Bollob\'{a}s-Riordan polynomial]
{Knot invariants and the Bollob\'{a}s-Riordan polynomial of embedded graphs}

\author{Iain Moffatt}
\address{Department of Applied Mathematics  \\
Charles University \\
Malostransk\' e n. 25 \\
118 00 \mbox{Praha 1} \\
Czech Republic.}

\email{iain@kam.mff.cuni.cz}

\thanks{
\newline
2000 {\em Mathematics Classification.} 57M15.
}

\date{\today}

\begin{abstract}
For a graph $G$ embedded in an orientable surface $\Sigma$, we consider associated links $\calL(G)$ in the thickened surface $\Sigma \times I$. We relate the HOMFLY polynomial of $\calL(G)$ to the recently defined \BR polynomial of a ribbon graph. This generalizes celebrated results of Jaeger and Traldi.
We use knot theory to prove results about graph polynomials
and, after discussing questions of equivalence of the polynomials, we  go on to use our formulae to prove a duality relation for the \BR polynomial. We then consider the specialization to the Jones polynomial and recent results of Chmutov and Pak to relate the \BR polynomials of an embedded graph and its tensor product with a cycle.
\end{abstract}

\maketitle


\section{Introduction}\label{sec:1}

The use of combinatorial methods in knot theory, as is well known, has  led to connections between graph and link polynomials, the most famous of which is a relation between the Jones polynomial and the Tutte polynomial. This appeared early on in the development of quantum topology (see \cite{kauf,this,jones}). Jaeger also found a connection between the Tutte and HOMFLY polynomials (\cite{Ja}). The relation was extended to a larger class of links by Traldi in \cite{Tr}.
 Here we are interested in the extension of these relationships to polynomials of embedded graphs and
 invariants of links in $\Sigma \times I$, the product of a surface $\Sigma$ and the unit interval $I = [0,1]$. We refer the reader to  \cite{CRR,hp,ik,Li,tu} for various approaches to the generalization of invariants of links in  the $3$-sphere $S^3$ to  links in $\Sigma \times I$ and to \cite{aflt,CP2,Da,dk,m} for recent related results.
In the recent preprint \cite{CP}  (subsequently revised to become paper \cite{CP2}) Chmutov and Pak  generalized the connection between the Jones and Tutte polynomials by relating the \BR polynomial of a ribbon graph $F$, which generalizes the Tutte polynomial to embedded graphs, to the Jones polynomial of a link in $F\times I$.

Rather than using graph theory to develop knot theory, as would be more usual, in this note we are interested in doing things the other way round: we use knot theory to advance our knowledge of graph polynomials. Our approach is to relate  the \BR polynomial of an embedded graph $G \subset \Sigma$ to the HOMFLY polynomial of a link in $\Sigma \times I$ and then to use the topology of knots to deduce results for the graph polynomial. This relation between polynomials  generalizes the results of Jaeger and Traldi mentioned above and also answers a question of Chmutov and Pak posed in \cite{CP}.
 We will show that the HOMFLY and the \BR polynomials are equivalent along the surface $\{ (x,y,z)|xyz^2=1 \}$. This observation allows us to answer graph theoretical questions with knot theory. We use  knot theory to prove a duality relation for the \BR polynomial (see also \cite{EMS}). Finding such a duality relation was a problem posed by Bollob\'{a}s and Riordan in \cite{BR}.
We also give an application to knot theory by showing that the genus of the smallest surface containing a link projection can be recovered from its HOMFLY polynomial.
We then go on to consider the  Jones polynomial as a specialization of the HOMFLY, and the relation of this with  Chmutov and Pak's result in \cite{CP}. From this we find a relation between the \BR polynomial of an embedded graph and its tensor product with the $(2^p+1)$-cycle $C_{2^p+1}$ (see also \cite{JVW,Wo}).

In section~\ref{sec:2} we recall how embedded graphs give rise to ribbon graphs and we define the \BR polynomial. We then discuss how to construct oriented links in $\Sigma \times I$ from an embedded graph $G \subset \Sigma$.
In section~\ref{sec:3}, after defining the HOMFLY polynomial of a link in $\Sigma \times I$, we prove our first result which relates it to the \BR polynomial.
This is generalized in two directions in section~\ref{sec:4} by considering weighted graphs and a multivariate \BR polynomial. We then go on to discuss the question of equivalence between the two polynomials and  prove a duality relation for the \BR polynomial.
The final section concerns the specialization to the Jones polynomial. We relate our earlier formulae for the HOMFLY polynomial to Chmutov and Pak's results on the Jones polynomial and prove a formula for the tensor product of a graph and its \BR polynomial.


\section{Links and embedded graphs}\label{sec:2}

\subsection{Ribbon graphs}
\begin{figure} \begin{center}
 \epsfig{file=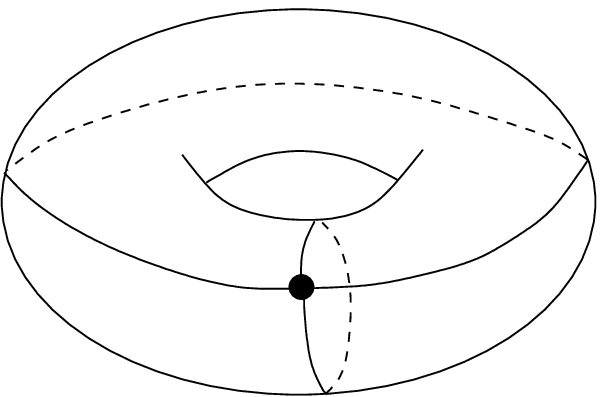, height=2cm}
\hspace{1.5cm}
\epsfig{file=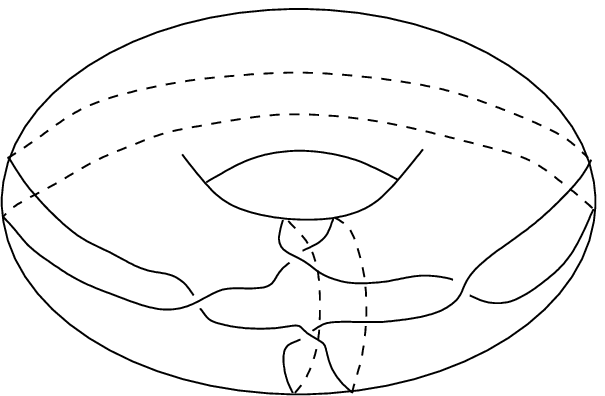, height=2cm}
\hspace{1.5cm}
\epsfig{file=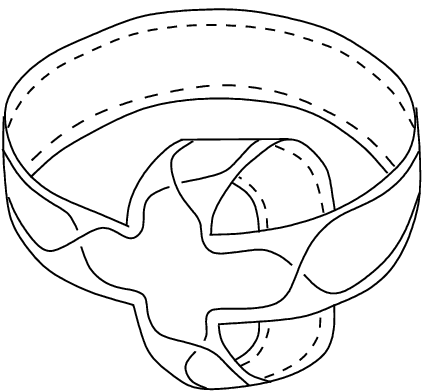, height=2cm} \end{center}
\caption{}
\label{fig:links}
\end{figure}

A {\em ribbon graph} (or {\em fatgraph})  $F$ is a graph together with a fixed cyclic order of the incident half-edges at each vertex. Ribbon graphs can be regarded as orientable surfaces with boundary by fattening  each vertex  into a disk, $D^2$, and fattening each edge into an `untwisted ribbon'.
Notice that if $G \subset \Sigma$ is a graph embedded in an orientable surface $\Sigma$, then a ribbon graph $F$ arises naturally as a neighbourhood of $G$ (we  retain information on the position of the vertices). This is indicated in figure~\ref{fig:links}. In such a situation we will say that $F$ is the ribbon graph {\em associated} with the embedded graph $G$. We will generally denote a ribbon graph by $F$ and an embedded graph by $G$ and  move freely between the two concepts.

\begin{figure}
\[ \epsfig{file=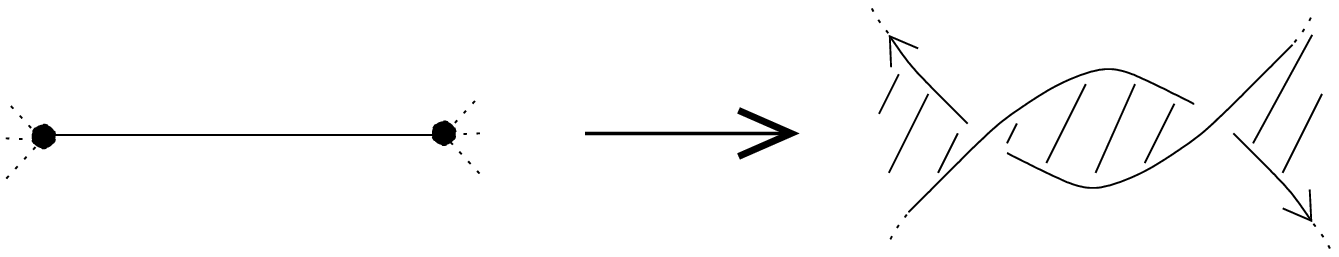, height=1.5cm} \]
\caption{}
\label{fig:jaeger}
\end{figure}

We are interested in the \BR polynomial of ribbon graphs (\cite{BR1,BR}) which is the natural generalization of the Tutte polynomial to embedded graphs.
We need some notation. Let $G=(V,E)$ be a graph then  we set $v(G)= |V|$,
$e(G)=|E|$, $r(G)= |V|- k(G)$ and $n(G)= |E|- r(G)$ where $k(G)$  denotes the number of connected components
of $G$. We  use a similar notation for ribbon graphs. In addition, if $F$ is a ribbon graph regarded as a surface, then we set $p(F)=|\partial(F)|$, the number of its boundary components. A {\em state} of a (ribbon) graph $F$ is a spanning sub(ribbon)~graph and we  denote the set of states by $\s(F)$. The {\em \BR polynomial} of a ribbon graph $F$ can then be defined as the sum over states:
\begin{equation}\label{eq:BRpoly}
R(F;\al , \be, \ga)  = \sum_{H\in \s(F)}\al^{r(F)-r(H)}\be^{n(H)}\ga^{k(H)-p(H)+n(H)}.
\end{equation}
By the \mbox{\BR} polynomial of an embedded graph we mean the \BR polynomial of the associated ribbon graph.

Notice that
$R(F; \al-1, \be-1 , 1)= T(F; \al , \be)$  the Tutte polynomial,
and, since the exponent of $\ga$ is exactly twice the genus of the ribbon graph regarded as a surface (see \cite{BR,LM}), 
$R(G; \al-1, \be-1 , \ga)= T(G; \al , \be)$ whenever $G \subset \mathbb{R}^2$.

We will also make use of the following rearrangement of the \BR polynomial which is obtained by expanding the rank  and nullity
\begin{equation}\label{eq:br}
 R(F;\al , \be, \ga)= \al^{-k(F)} (\be \ga)^{-v(F)} \sum_{H\in \s(F)} (\al\be\ga^2)^{k(H)} (\be\ga)^{e(H)} (\ga)^{-p(H)}.
\end{equation}

\subsection{Links in $\Sigma \times I$}
Having discussed embedded graphs,  we move onto our second main object.
A {\em link} is an embedding of a finite number of copies of the unit circle $S^1$ into a 3-manifold $M$. Knot theory is mostly concerned with the special case $M=S^3$. Here, however, we are concerned with the more general case of links in $\Sigma \times I$, where $\Sigma$ is an orientable surface and $I=[0,1]$ is the unit interval. Note that knot theory in the $3$-ball $D^2 \times I$ is equivalent to knot theory in the $3$-sphere $S^3$.
Given a link $L \subset \Sigma \times I$ there is a generic projection on to the surface $\Sigma$ by projection onto the first variable. It is obvious how to construct  a link in $\Sigma \times I$  from its projection.

If we are given an embedded graph $G \subset \Sigma$  we can associate an oriented link projection by replacing each edge with the oriented tangle of figure~\ref{fig:jaeger} and connecting these tangles according to the cyclic order at the vertices. An example is shown in figure~\ref{fig:links}. We denote a link in $\Sigma \times I$ constructed in this way by $\calL(G)$.
We may assume that the link lies in a sufficiently small neighbourhood of $G$ so that we obtain a projection $\calL(F)$ on the associated  ribbon graph $F$. This gives rise to a link in $F\times I$. Again this is indicated in figure~\ref{fig:links}.


\section{The HOMFLY in $F \times I$}\label{sec:3}

\subsection{The HOMFLY polynomial}
\begin{figure}
\[ \begin{array}{ccccc}
\epsfig{file=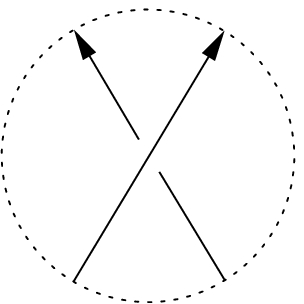, height=1.5cm} & & \epsfig{file=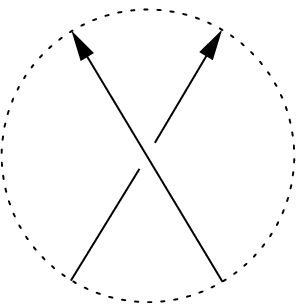, height=1.5cm} & & \epsfig{file=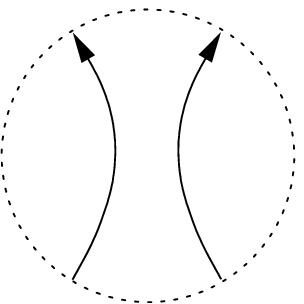, height=1.5cm} \\
L_+ & & L_- & & L_0
\end{array}\]
\caption{}
\label{fig:skein}
\end{figure}

Let $\Sigma$ be an orientable surface (possibly with boundary) and $L\subset \Sigma \times I$ be a link.
The HOMFLY polynomial $P(L)$ is a link invariant which satisfies the {\em HOMFLY skein relation}
\[ xP(L_+) -x^{-1}P(L_-) = yP(L_0) ,\]
where $L_+$,  $L_-$ and $L_0$ are links which differ only in the locality of a single crossing as shown in figure~\ref{fig:skein}.
We will also require that the invariant takes the value $1$ on the unknot, {\em i.e.} $ P(\mathcal{O})=1$.

If we were only considering links in $S^3$ or $D^2 \times I$ then this would be enough to uniquely determine a link invariant $P(L)\in \Z [x^{\pm 1}, y^{\pm 1}]$.
 However, for a general orientable surface $\Sigma$ this is not enough. To describe a basis for the HOMFLY skein module we need to introduce the notion of a descending link.

First notice that there is a natural {\em product} of links in $\Sigma\times I$ given by reparameterizing the two copies of  $\Sigma\times I$  and stacking them:
\[
 (\Sigma \times I) \times  (\Sigma \times I) \cong (\Sigma \times [1/2,1]) \times  (\Sigma \times [0,1/2])
\rightarrow (\Sigma \times I) .
\]
We denote the projections from $\Sigma \times I$ to $\Sigma $ and $I$ by $p_\Sigma $ and $p_I$ respectively.  The value $p_I(x)$ is called the {\em height} of $x$. We can now make our definition.
\begin{definition}

(1) A knot $K \subset \Sigma  \times I$ is {\em descending} if
it is isotopic to a knot $K^{\prime} \subset \Sigma  \times I$ with the property that
 there is a choice of basepoint $a$ on $K^{\prime}$ such that if we travel along $K^{\prime}$ in the direction of the orientation from the basepoint the height of $K^{\prime}$ decreases until we reach a point $a^{\prime}$ with $p_\Sigma (a)=p_\Sigma (a^{\prime})$ from which $K^{\prime}$ leads back to $a$ by increasing the height and keeping the projection onto $F$ constant.

(2) A link $L \subset \Sigma  \times I$ is said to be {\em descending} if it is isotopic to a product of descending knots.
\end{definition}
Clearly each (isotopy class of a) descending knot uniquely determines  a conjugacy class of the fundamental group $\pi_1(\Sigma )$.
Moreover there is a bijection between the conjugacy classes in $\pi_1(\Sigma )$ and isotopy classes of descending knots. In other words a conjugacy class determines a descending knot. 
 In \cite{Li}, Lieberum gives a procedure for choosing a canonical element of the conjugacy classes of $\pi_1(\Sigma )$. We will denote this set of choices by $S(\Sigma ) = \{t_w\}$. We do not need to worry about the exact choices here.
This means that each descending knot $K$ determines some $t_K \in S(\Sigma)$ and each descending link determines a monomial $t_L$ in commuting indeterminates $\{t_w\}$ (we have $t_{\mathcal{O}}=1$).
For a descending link $L$ we set
\begin{equation}\label{eq:basis}
P(L) = t_L \left( \frac{x-x^{-1}}{y} \right)^{k(L)-1},
\end{equation}
where $k(L)$ is the number of components of the link $L$.

The {\em HOMFLY} (or {\em HOMFLYPT}) polynomial is then defined by the following theorem.
\begin{theorem}[\cite{Li}]
There exists a unique  invariant $P(L) \in \Z [x^{\pm 1}, y^{\pm 1}]\otimes \Z [\{t_w\}]$ of links $L \subset \Sigma \times I$ that satisfies the HOMFLY skein relation and equation~\ref{eq:basis}.
\end{theorem}
If we set each $t_w=1$ then we obtain a polynomial in $\Z [x^{\pm 1}, y^{\pm 1}]$, which we denote by $P(L; x,y )$.
We consider this $2$-variable polynomial first.


\subsection{The relation to the Bollob\'{a}s-Riordan polynomial}
In our first result we relate the \BR polynomial of an embedded graph $G \subset \Sigma$ to the HOMFLY polynomial of the associated link $\calL(G) \subset \Sigma \times I$. In order to do this we need to set each of the variables $t_w$ equal to $1$ and consider the invariant $P(L; x,y )$. In the following section we will discuss the full HOMFLY invariant $P(L) \in \Z [x^{\pm 1}, y^{\pm 1}]\otimes \Z [\{t_w\}]$ and relate it to a multivariate \BR polynomial. The relation proved in this section will follow from the more general result, however we feel that it is more clear if we give a separate proof.

As promised, the following theorem gives a relationship between the HOMFLY and \BR polynomials.
\begin{theorem}\label{th:homfly}
Let $\Sigma$ be an orientable surface, possibly with boundary, and $G \subset \Sigma$ be an embedded graph. Then
\begin{multline}\label{eq:homfly}
P(\calL(G); x,y ) = \left( \frac{1}{xy} \right)^{v(G)-1}
 \left( \frac{y}{x}\right)^{e(G)}
\left(  x^{2}-1  \right)^{k(G)-1}
 \\ R\left(G; x^{2}-1 , \; \frac{x-x^{-1}}{xy^2},\; \frac{y}{x - x^{-1}}\right).
\end{multline}
\end{theorem}
\begin{proof}
An application of the HOMFLY skein and isotopy of the link gives the equation
\begin{equation}\label{eq:skstate}
P\left( \raisebox{-5mm}{\epsfig{file=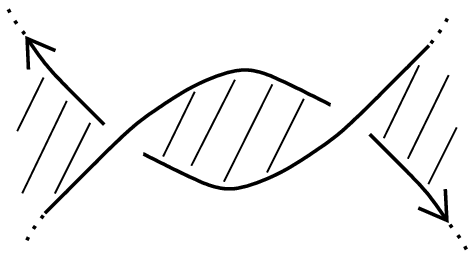 ,height=1cm}}  \right) =
\frac{1}{x^2} \; P\left( \raisebox{-5mm}{\epsfig{file=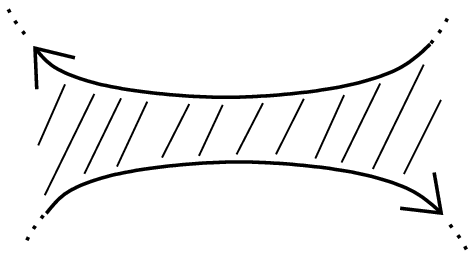 ,height=1cm}}  \right)
+\frac{y}{x} \; P\left( \raisebox{-5mm}{\epsfig{file=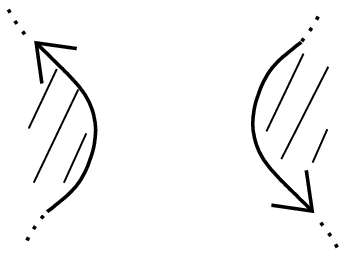 ,height=1cm}}  \right).
\end{equation}
We will call a link obtained by eliminating the crossings in this way a {\em resolution}.
There is a clear correspondence between resolutions of the link $\calL (G)$ and the states of the associated ribbon graph $F$ given by including an edge in the state whenever the corresponding link is resolved as
\raisebox{-2mm}{\epsfig{file=figs/st1 , height=0.6cm}}
and excluding an edge otherwise.

The links corresponding to a given state have no crossings and are therefore descending. It is then clear that
\begin{equation*}\begin{split}
P(\calL(G); x,y ) & =
 \sum_{H \in \s (G)}  \left( \frac{1}{x^2} \right)^{e(H)}
\left( \frac{y}{x} \right)^{e(G)-e(H)}
\left( \frac{x-x^{-1}}{y} \right)^{p(H)-1}\\
& =\left( \frac{y}{x-x^{-1}} \right) \left( \frac{y}{x} \right)^{e(G)}
\sum_{H \in \s (F)}  \left( \frac{1}{xy} \right)^{e(H)}\left( \frac{x-x^{-1}}{y} \right)^{p(H)}.
\end{split}
\end{equation*}
Now setting $\al = x^2-1$, $\be=(x-x^{-1})/xy^2$ and $\ga=y/(x-x^{-1})$, so that $\al\be\ga^2 = 1$ and $\be\ga=1/xy$, we see that by applying (\ref{eq:br}), the above can be written as
\begin{equation*}\begin{split}
\left(\frac{y}{x-x^{-1}}\right)
\left( \frac{y}{x} \right)^{e(G)}
(x^2-1)^{k(G)}
\left( \frac{1}{xy} \right)^{v(G)}
\; R\left(F; x^{2}-1 , \; \frac{x-x^{-1}}{xy^2},\; \frac{y}{x - x^{-1}}\right)
\\
=\left( \frac{1}{xy} \right)^{v(G)-1}
 \left( \frac{y}{x}\right)^{e(G)}
\left(  x^2 -1  \right)^{k(G)-1}
 \; R\left(F; x^{2}-1 , \; \frac{x-x^{-1}}{xy^2},\; \frac{y}{x - x^{-1}}\right)
\end{split}
\end{equation*}
as required. 
\end{proof}
Notice that when $\Sigma$ is a disc then, using the relation between the \BR and   Tutte polynomials, we recover Jaeger's expression for the HOMFLY from \cite{Ja}.


\section{The full polynomial}\label{sec:4}
In order to find a formula for the full HOMFLY invariant, rather than just its specialization at $t_w=1$, we need to consider a ``Sokalization'' ({\em c.f.} \cite{So}) of the \BR polynomial.
The \BR polynomial can be expressed as
\begin{equation*}\begin{split}
R(F;\al,\be,\ga)
& = \al^{-k(F)} (\be \ga)^{-v(F)} \sum_{H \in \s (F)} (\al\be\ga^2)^{k(H)} (\be\ga)^{e(H)} (\ga^{-1})^{p(H)} \\
& = ((ac)^{-1}b)^{k(F)} b^{-v(F)} \sum_{H \in \s (F)} a^{k(H)} b^{e(H)} c^{p(H)},
\end{split}
\end{equation*}
where $a=\al\be\ga^2$, $b=\be\ga$ and $c=\ga^{-1}$. Therefore the \BR polynomial can be equivalently formulated as a polynomial
\[
B(F; a,b,c) = \sum_{H \in \s (F)} a^{k(H)} b^{e(H)} c^{p(H)}.
\]

A {\em weighted ribbon graph} is a ribbon graphs equipped with a map  from its edge set $E$ to a set of {\em weights} $\{b_e \}_{e \in E(F)}$. Notice that we can (and will) regard the {\em weights} $b_e$ as formal commuting variables. We can regard an unweighted ribbon graph as a weighted graph all of whose weights are equal. We will usually denote this single weight by $b$.

For a weighted ribbon graph $F$ we define the {\em weighted \BR polynomial} as
\[
B(F; a,\vb ,c) = \sum_{H \in \s (F)} a^{k(H)} \left(\prod_{e\in E(H)} b_e \right) c^{p(H)},
\]
where $\vb = \{b_e \}_{e \in E}$ is the set of edge weights. We will sometimes exclude the variables from the notation and just write $B(F)$. As before, for an embedded graph $G$, $B(G):=B(F)$ where $F$ is the associated ribbon graph.
Notice that for an unweighted ribbon graph this polynomial is equivalent to  the \BR polynomial.

We now state the general form of Theorem~\ref{th:homfly}.
\begin{theorem}\label{th:fullhomfly}
Let $\Sigma$ be an orientable surface, possibly with boundary, and $G \subset \Sigma$ be an embedded graph
 with edge weights $\vb = \{b_e \}_{e \in E(F)} $. Then
\begin{multline}\label{eq:full} 
P(\calL(G); x,y, \{ t_w \} ) =\\  f \left( \
\begin{array}{c}
\left( y/(x-x^{-1}) \right)
 \left( y/x\right)^{e(G)}
 \; B\left(G; 1 , \; \left(b_e/xy\right)_{e\in E(G)},\; (x - x^{-1})/y\right)
 \end{array}
 \right)
\end{multline}
where $f$ is a function from polynomials in $\vb$ to polynomials in the conjugacy classes $S(\Sigma) =\{t_w\}$ of $\pi_1 (\Sigma)$.
\end{theorem}

\begin{proof}
As in the proof of Theorem~\ref{th:homfly}, there is a correspondence between the states of the ribbon graph and the resolutions of the link shown in Equation~\ref{eq:skstate}. So it is enough to show that each corresponding term on the left and right hand sides of (\ref{eq:full}) is assigned the same value.

Let $F$ be the ribbon graph associated with $G$. Observe that a  state $H$ of the ribbon graph $F$ is assigned the monomial in $\vb$
\[
 \left( \frac{y}{x-x^{-1}} \right)
 \left( \frac{y}{x}\right)^{e(F)}
  \left(\prod_{e\in E(H)} \frac{b_e}{xy} \right)
\left(\frac{x - x^{-1}}{y} \right)^{p(H)} \]
by the right hand side. This is equal to
\begin{equation}\label{eq:fullpf1}
\left( \frac{y}{x}\right)^{e(F)-e(H)}
\left( \frac{1}{x^2}  \right)^{e(H)}
\left(\frac{x - x^{-1}}{y} \right)^{p(H)-1}
 \left(\prod_{e\in E(H)} b_e \right).
\end{equation}

Now $H$ is a state of $F$ so we can regard it as a subsurface of $F$. The boundary $\partial H$ determines a descending link $L_H \subset F\times I$.
Recalling the correspondence between states of $F$ and resolutions of $\calL (F)$ in the proof of Theorem~\ref{th:homfly}, it is clear that $L_H$ is a resolution of $\calL (F)$. By (\ref{eq:skstate}) and (\ref{eq:basis}) this resolution contributes the expression
\[
\left( y/x\right)^{e(F)-e(H)}
\left( 1/x^2  \right)^{e(H)}
\left((x - x^{-1})/y \right)^{p(H)-1}
 t_{L_H}.
\]
Comparing this to (\ref{eq:fullpf1}) we see that all that remains is to describe the map $f: \prod_{e\in E(H)} b_e  \mapsto t_{L_H}$ which induces the map of the theorem.
Now the monomial $\prod_{e\in E(H)} b_e$ tells us which edges are in the state $H$ and therefore uniquely determines the subsurface $H$
 of $F$ and the link diagram $L_H$ on the surface $\Sigma$. 
We simply define $f$ to be the map which assigns the appropriate conjugacy class of $\pi_1 (\Sigma)$ and if a monomial contains an element $b_e$ more than once then send it to zero. 
\end{proof}

\begin{remark}
It is a simple exercise to give an explicit construction of the map $f$ into the representatives of the conjugacy classes described in \cite{Li}. All that is required is the observation that any ribbon graph $F$ is homeomorphic to a decomposed surface of \cite{Li} then simply use this homeomorphism to pull back the link $L_H$ to the set of canonical conjugacy generators of the decomposed surface.
\end{remark}


\subsection{On Traldi's extension}

So far we have discussed the construction of links using Jaeger's idea of replacing an edge of the ribbon graph as in figure~\ref{fig:jaeger}. However in \cite{Tr} Traldi  extended this idea by replacing the edges of a weighted graph by various tangles according to the weight. We use this idea to extend Theorem~\ref{th:fullhomfly}.

\begin{figure}
\[\begin{array}{ccccccc}
\epsfig{file=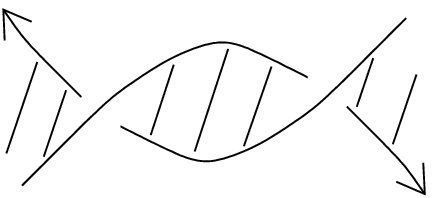, height=1.2cm} &&
\epsfig{file=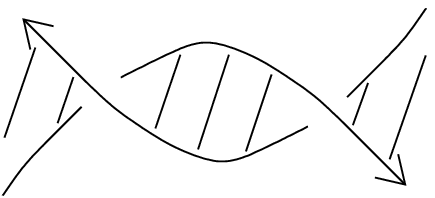, height=1.2cm} &&
\epsfig{file=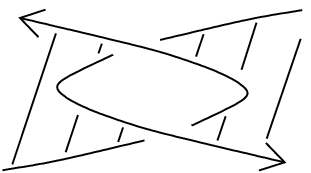, height=1.2cm} &&
\epsfig{file=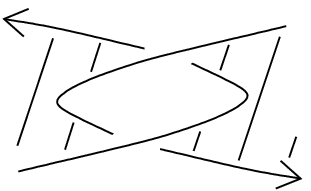, height=1.2cm}
\\
w_1 & &w_2 & &w_3 & &w_4
\end{array}\]
\caption{}
\label{fig:traldi}
\end{figure}

Let $F$ be a ribbon graph with edge weights in the set
$\{ b_e \}_{e\in E(F)}\times \{ w_1, \ldots , w_4 \}$.
Construct a link $\calL (F)$ by associating a tangle to each edge according to the weight $w_i$ as indicated in figure~\ref{fig:traldi}.
The following result then holds.
\begin{theorem}
Let $\Sigma$ be an orientable surface, possibly with boundary, and $G \subset \Sigma$ be an embedded graph
 with edge weights
 $\{ b_e \}_{e\in E(F)}\times \{ w_1, \ldots , w_4 \}$. Then
\begin{multline}\label{eq:traldi}
P(\calL(G)) = f \left( 
\begin{array}{c}
 \left( y/(x-x^{-1}) \right)
 \left( y/x \right)^{w_1(G)}
\left( -yx \right)^{w_2(G)}
\left( 1/x^2 \right)^{w_3(G)}
\left( x^2 \right)^{w_4(G)}
\end{array}
\right. \\
  \left. 
  \begin{array}{c}
  B\left(G; 1 , \; \left( b_e  w_e \right)_{e\in E(G)},\; (x - x^{-1})/y\right)
  \end{array}
  \right)
\end{multline}
where $w_1 = 1/xy$, $w_2=-x/y$, $w_3=xy$ and $w_4=-y/x$;
$f$ is a function from polynomials in $\vb$ to polynomials in  the conjugacy classes $S(\Sigma) =\{t_w\}$ of $\pi_1 (\Sigma)$ and $w_i(G)$ denotes the number of edges of $G$ with the edge weight of the form $(\cdot , w_i)$.
\end{theorem}
The proof of this result is a straight forward generalization of the proof of Theorem~\ref{th:fullhomfly} and is therefore excluded.

Notice that Traldi's Theorem~5 in \cite{Tr} and also Theorem~\ref{th:fullhomfly} (by setting all $w_e = w_1$) may be recovered from this.


\subsection{Determination of the \BR polynomial and a duality relation}

For planar graphs, Jaeger showed that  the Tutte polynomial of a graph $G$ and the HOMFLY polynomial of the associated link $\calL (G)$ are equivalent on domain $\{ (u,v)\in \mathbb{R}^2 | u \neq 0 \}$. Since the \BR polynomial is a generalization of  the Tutte polynomial for non-planar embedded graphs, it is natural to ask if the HOMFLY polynomial of $\calL (F)$ determines the \BR polynomial of $F$. We will see that this isn't quite the case.

To answer this question we need to find non-zero values of  $x$ and $y$ such that
$\al = x^2-1$, $\be = (x-x^{-1})/xy^2$ and $\ga = y/(x-x^{-1})$. We see that this is possible for the choices $x= \sqrt{\al+1} \neq 1$ and
 $y= \sqrt{\al/(\be(\al +1))}$. Notice that $\al \be \ga^2 =1$. 
 The proof of Proposition~\ref{pr:det} then follows immediately.
\begin{proposition}\label{pr:det} For any ribbon graph $F$ and $\al \neq 0,1$, $\be \neq 0$
\begin{multline}\label{pr:determination}
R\left( F; \al, \be , 1/\sqrt{\al \be} \right) = 
 \left(  \sqrt{\al}/ \sqrt{\be}  \right)^{v(F)-e(F)-1}
(\al+1)^{e(F)} \al^{1-k(F)} \\
P\left(\calL(F);  \sqrt{\al-1} , \sqrt{\frac{\al}{\be (\al+1)}} \right).
\end{multline}
\end{proposition}

We observe the following intriguing application of the above to knot theory. 
We say that a link $L \subset \Sigma \times I$ is {\em essential} in $\Sigma$ if there is no embedded surface $\Sigma^{\prime} \subset \Sigma$, whose genus is smaller than that of $\Sigma$, such that $L \subset \Sigma^{\prime}\times I$.
Also we say that a link $L \subset \Sigma \times I$ is {\em split} if it is isotopic to a link $L^{\prime}$ with the property that $ L^{\prime} \cap (\Sigma \times [0,1/2]) \neq \emptyset$, $L^{\prime} \cap (\Sigma \times [1/2,1])\neq \emptyset$ and $L^{\prime}  \cap (\Sigma \times \{  1/2 \}) = \emptyset$.
\begin{corollary}
Let $L \subset \Sigma \times I$ be a non-split alternating link essential in $\Sigma$. If  $L$ has a projection of the form $\calL(G)$ for some graph $G \subset \Sigma$, then the genus of $\Sigma$ can be recovered from the HOMFLY polynomial of $L$.
\end{corollary}
\begin{proof}
Since $L$ has a projection of the form $\calL(G)$, $G \subset \Sigma$, and the HOMFLY polynomial is equivalent to the \BR polynomial on the surface $\{ (x,y,z)|xyz^2=1 \}$, we see it is enough to show that the genus of $\Sigma$ can be recovered from $R\left( G; \al , \be, 1/\sqrt{\al\be}  \right)$.

Now by \cite{LM},
\[ R\left( G; \al , \be, 1/\sqrt{\al\be}  \right)
= \sum_{H\in \s(F)}\al^{r(F)-r(H)}\be^{n(H)} (\al \be)^{-g(H)}.
\]
By setting $\be=1$ and observing that $G$ must be connected  the above sum becomes   \\ $\sum_{H\in \s(G)} \al^{k(H)-1 -g(H)}$. Therefore the lowest degree of $\al$ is equal to  $-g(\Sigma)$ (since $L$ is essential in $\Sigma$), and we are done. 
\end{proof}
This corollary motivates the following conjecture.
\begin{conjecture}
Let $L \subset \Sigma \times I$ be a link. Then the genus of an essential surface for $L$ can be recovered from its HOMFLY polynomial $P(L)$.
\end{conjecture}

\begin{remark}
We consider the 2-variable HOMFLY skein. However there is also a 3-variable version of the skein relation. It seems reasonable to conjecture that the 3-variable HOMFLY polynomial has the basis
$ t_L \left( -(x+y)/z \right)^{k(L)-1}$ on descending links. In which case one could make similar arguments as above and relate the 3-variable HOMFLY and \BR polynomials. However, with regards  to the proposition above, the additional variable in the HOMFLY does not provide any further information about the \BR polynomial. In fact all that happens is that we would introduce a redundant third variable into the right hand side of equation~\ref{pr:determination}. This is similar   to  Proposition~2 of \cite{Ja}.
\end{remark}

\medskip

Given a graph $G$ embedded in a 2-manifold without boundary, one can form a {\em dual} embedded graph $G^*$ in the usual way. From this we can obtain the dual $F^*$ of a ribbon graph $F$.  In the remainder of this section we study the relation between the \BR polynomial of a ribbon graph and its dual.

As far as the author is aware, there is no known duality relation for the full 3-variable \BR polynomial. In \cite{BR}, Bollob\'{a}s and Riordan prove a 1-variable relation
$R(F; \al, \al ,\al^{-1}) =  R(F^*; \al, \al ,\al^{-1})$ leaving it as an open problem to find a multi-variable relation. Ellis-Monaghan and Sarmiento in \cite{EMS} extended this to the  2-variable relation:
\begin{theorem}[\cite{EMS}] \label{th:duality}
Let $F$ be a connected ribbon graph and $F^*$ its dual. Then
\[
R(F; \al, \be , 1/\sqrt{\al \be}) = \left( \be /\al\right)^{g(F)} R(F^*; \be, \al , 1/\sqrt{\al \be}),\]
where $g(F)$ is the genus of the ribbon graph regarded as a surface.

\end{theorem}
One notices immediately that the specialization of the \BR polynomial  in this relation is exactly that which is determined by the HOMFLY polynomial in (\ref{pr:determination}).  We will provide a new proof for the duality relation above using knot theory and we will see that the duality relation holds for the specialization $R(F; \al, \be , 1/\sqrt{\al \be})$ precisely because it is determined by the HOMFLY polynomial.

\begin{proof}[Proof of Theorem~\ref{th:duality}.]
Let $\calL_1 (F)$ be the link associated the ribbon graph $F$ by associating tangles $w_1$ of figure~\ref{fig:traldi} to edges of $F$ and let $\calL_3 (F^*)$ be the link associated the ribbon graph $F^*$ by associating tangles $w_3$ of figure ~\ref{fig:traldi} to edges  and reversing the orientation of all components of the link. Clearly these two links are isotopic and therefore
$P(\calL_1 (F))=P(\calL_3 (F^*))$.
Now by (\ref{eq:homfly}) we have
\begin{equation*}
P(\calL_1(F); x,y ) = \left( xy \right)^{v(G)-1}
 \left( 1/x^2\right)^{e(G)}
 \; R\left(F; x^{2}-1 , \; \frac{x-x^{-1}}{xy^2},\; \frac{y}{x - x^{-1}}\right),
\end{equation*}
and since the reversal of the orientation of a link does not change its HOMFLY polynomial, equation~\ref{eq:traldi} gives
\begin{equation*}
P(\calL_3(F^*); x,y ) = \left( \frac{1}{xy} \right)^{v(G)-1}
 \left( \frac{y}{x}\right)^{e(G)}
 \; R\left(F^*; \; \frac{x-x^{-1}}{xy^2}, x^{2}-1 , \; \frac{y}{x - x^{-1}}\right).
\end{equation*}
Then, by the isotopy of the link, we have
\begin{multline*}
R\left(F; x^{2}-1 , \; \frac{x-x^{-1}}{xy^2},\; \frac{y}{x - x^{-1}}\right) = \\
\left( xy \right)^{v(F)+v(F^*)-e(F)-2}
R\left(F^*; \; \frac{x-x^{-1}}{xy^2}, x^{2}-1 , \; \frac{y}{x - x^{-1}}\right).
\end{multline*}
Notice that since $v(F^*)$ equals the number of faces of $F^*$ (embedded in a non-punctured surface), we have $ v(F)+v(F^*)-e(F)-2 = \chi (F) -2 =-2g(F)$.
Finally the substitutions
$x= \sqrt{\al+1} \neq 1$ and $y= \sqrt{\al/(\be(\al +1))}$ give the relation
\[
R(F; \al, \be , 1/\sqrt{\al \be}) = \left( \be /\al\right)^{g(F)} R(F^*; \be, \al , 1/\sqrt{\al \be})
\]
as required. 
\end{proof}


\section{The Jones polynomial}\label{sec:5}

Given a digraph $G$, a graph $H$ and a distinguished oriented edge of $H$. The {\em tensor product} $G \otimes H$ is defined to be the graph obtained by identifying each edge of $G$ with the distinguished edge of a copy of $H$ and then deleting each of the edges of $G$ ({\em i.e.} we take the 2-sum with $H$ at every edge of $G$). In general the graph obtained depends upon the various choices made. However this is not always the case. One example for which the tensor product is well defined is when $H$ is the $p$-cycle $C_p$. In this case notice also that the tensor product is independent of the orientation of the edges of $G$ and so the tensor product makes sense for (embedded) graphs $G$.
Tensor products and their effect on the Tutte polynomial have been considered previously in \cite{JVW} and \cite{Wo}.
In this section we consider the connection between the tensor product $G\otimes C_3$, Theorem~\ref{th:homfly} and a result of Chmutov and Pak which relates the Jones polynomial of a link in $F\times I$ to the \BR polynomial. This generalizes results of Huggett which appeared in \cite{Hu}.

\begin{figure} \begin{center}
\epsfig{file=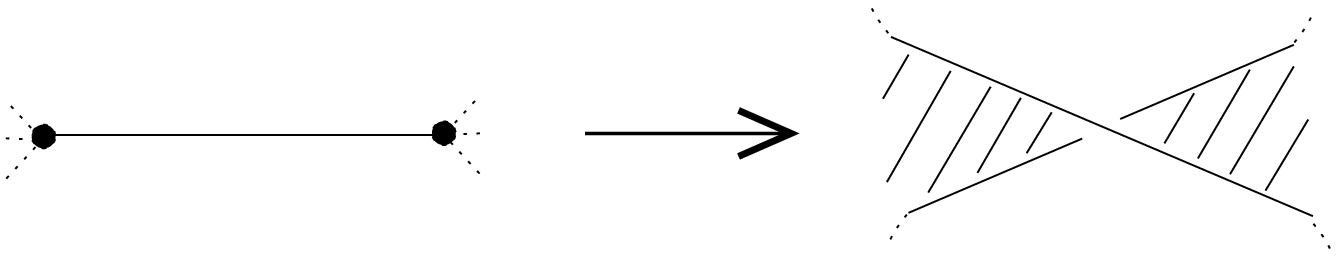, height=1.2cm} \end{center}
\caption{}
\label{fig:medial}
\end{figure}

So far we have constructed oriented links from embedded graphs by replacing each edge with a tangle as in figures~\ref{fig:jaeger}~and~\ref{fig:traldi}. However an unoriented link can be associated to an embedded graph by replacing each edge with the tangle indicated in figure~\ref{fig:medial}. This is known as the {\em medial link} and we  denote the medial link associated with a graph $G$ by $L(G)$.

The Jones polynomial $J(L)$ of a link $L$ is defined by the skein
$t^{-1} J(L_+) - t J(L_-) = (t^{1/2} - t^{-1/2})J(L_0)$ and $J(\mathcal{O})=1$.
There is a well known formula relating the Tutte polynomial $T(G;-t, -t^{-1})$  of a planar graph $G$ to the Jones polynomial of its medial link. This result was recently generalized by Chmutov and Pak who related the \BR polynomial of a ribbon graph $F$ and the Jones polynomial of its medial link $L(F) \subset F\times I$:
\begin{theorem}[\cite{CP}]
Let $F$ be a ribbon graph and $L(F) \subset F\times I$ be its medial link then
\begin{multline}\label{eq:cp}
J(L(F); t)= (-1)^{\omega}t^{(3\omega - r(F) +n(F))/4} (-t^{1/2} - t^{-1/2})^{k(F)-1} \\
R\left(F; -t-1 , -t^{-1}-1, 1/(-t^{1/2} - t^{-1/2})\right),
\end{multline}
where $\omega$ is the writhe of the link (the writhe is the number of $L_+$ crossings minus the number of $L_-$ crossings, where $L_{\pm}$ are as in figure~\ref{fig:skein}.
\end{theorem}
Chmutov and Pak proved this result by considering the Kauffman bracket construction of the Jones polynomial.

Of course the Jones polynomial is the specialization of the HOMFLY polynomial at $x=t^{-1}$ and $y=t^{1/2}-t^{-1/2}$ and so by (\ref{eq:homfly}) we also have
\begin{multline}\label{eq:joho}
 J(\mathcal{L}(F);t) =
(t^{1/2}-t^{-1/2})^{(e(F)-v(F)+1  )} t^{(e(F)+v(F)-1)}
(t^{-2}-1)^{k(F)-1} \\
\; R\left( F; \, t^{-2}-1 , \, \frac{1-t^{2}}{(t^{1/2}- t^{-1/2})^2}  , \,
\frac{1}{-t^{1/2}-t^{-1/2}}  \right).
 \end{multline}

Now tensors return to the story. It is clear that, forgetting the orientations, the links $L(F\otimes C_3)$ and the mirror image of $\calL(F)$ (obtained by reversing all of the crossings) are isotopic. Therefore $J(L(F\otimes C_3); t^{-1}) = J(\mathcal{L}(F);t)$ (since taking the mirror image of a link changes the Jones polynomial by the substitution $t \mapsto t^{-1}$), where the orientation on $L(F\otimes C_3)$ is induced by that of $\calL(F)$.  We use this observation to prove the following formula for the tensor product of a ribbon graph with $C_3$.
\begin{theorem}\label{th:tens}
 Let $F$ be a ribbon graph; $\al \neq -2, 0, 1$; $\be = \al (1-\al)$ and $\ga = 1/ \sqrt{\al\be}$. Then
\[R(F\otimes C_3;\alpha, \beta, \gamma) = (\alpha +2)^{n(F)} R(F; \alpha (\alpha +2) ,\beta/(\alpha +2) , \gamma).  \]
\end{theorem}
\begin{proof}
Consider the links $\calL (F)$ and  $L=L(F\otimes C_3)$ described above. We will use the relation between their Jones polynomials to prove the result.  For brevity we set $A:=F \otimes C_3$.

By (\ref{eq:cp}) we have
\begin{multline*}
J(L(A); t) =
(-1)^{\omega}t^{(3\omega(L) - r(A) +n(A))/4} (-t^{1/2} - t^{-1/2})^{k(A)-1} \\
R\left(A; -t-1 , -t^{-1}-1, 1/(-t^{1/2} - t^{-1/2})\right).
\end{multline*}
Clearly $\omega(L) = -e(A)=-2e(F)$, $e(A)=2e(F)$ and $v(A)=v(F)+e(F)$ and the above  becomes
\begin{equation}\label{eq:proof1}
t^{(-3e(F) -v(F) +k(F))/2} (-t^{1/2} - t^{-1/2})^{k(G)-1}
R\left(A; -t-1 , -t^{-1}-1, 1/(-t^{1/2} - t^{-1/2})\right).
\end{equation}
On the other hand equation~\ref{eq:joho}, coming from the HOMFLY polynomial gives
\begin{multline*} 
J(\calL(A); t^{-1})  = (t^{2}+1)^{k(F)-1}
t^{(-e(F) -v(F) +1)/2} (-t^{1/2}+t^{-1/2})^{e(F)-v(F)+1} \\
R\left( F; \, t^{2}-1 , \, \frac{1-t^{-2}}{(t^{-1/2}- t^{1/2})^2} -1 , \,
\frac{1}{-t^{1/2}-t^{-1/2}}  \right),
\end{multline*}
which can be written as
\begin{multline}\label{eq:proof2}
(t^{2}+1)^{k(F)-1}
t^{(-e(F) -v(F) +1)/2} (-t^{1/2}+t^{-1/2})^{e(F)-v(F)+1} \\
R\left( F; \, t^{2}-1 , \, \frac{1-t^{-2}}{(t^{-1/2}- t^{1/2})^2} -1 , \,
\frac{1}{-t^{1/2}-t^{-1/2}}  \right).
\end{multline}
As observed above, we have $J(L(F\otimes C_3); t) = J(\mathcal{L}(F);t^{-1})$ and therefore, using (\ref{eq:proof1}) and  (\ref{eq:proof2}), we have
\begin{multline*} 
R\left(A; -t-1 , -t^{-1}-1, 1/(-t^{1/2} - t^{-1/2})\right) = \\
(-t+1)^{n(F)}
R\left( F; \, t^{2}-1 , \, \frac{1-t^{-2}}{(t^{-1/2}- t^{1/2})^2} -1 , \,
\frac{1}{-t^{1/2}-t^{-1/2}}  \right).
\end{multline*}
Finally, substituting $t = 1-\al$ gives
\[R(A;\alpha, \beta, \gamma) = (\alpha +2)^{n(F)} R(F; \alpha (\alpha +2) ,\beta/(\alpha +2) , \gamma) , \]
where $\be$ and $\ga$ are as in the statement of the lemma.
\end{proof}
Induction gives the following corollary.
\begin{corollary}
Let $F$ be a ribbon graph; $\al \neq -2, 0, 1$; $\be = \al (1-\al)$ and $\ga = 1/ \sqrt{\al\be}$. Then
\[R(F\otimes C_{2^{p}+1};\alpha, \beta, \gamma) =  \left(\sum_{i=0}^{2^p} (\al +1)^i  \right)^{n(F)}
 R\left(F; \; (\alpha +1)^{2^p}-1 ,  \; \frac{\beta}{\sum_{i=0}^{2^p} (\al +1)^i } , \; \gamma \right).  \]
\end{corollary}

\begin{remark}
In fact, it is not difficult to prove that for any $\al$, $\be$ and $\ga$,
$R(F\otimes C_3;\alpha, \beta, \gamma) = (\alpha +2)^{n(F)} R(F; \alpha (\alpha +2) ,\beta/(\alpha +2) , \gamma)$. Taking this as our starting point, we see that formula~\ref{eq:cp}, which is proven by considering the Kauffman bracket, and formula~\ref{eq:joho} are related through the notion of the tensor product of an embedded graph. This generalizes the main result of \cite{Hu}.
\end{remark}


\begin{thebibliography}{[99]}
\bibitem{aflt} C. Adams, T. Fleming, TM. Levin and A. M. Turner,
Crossing number of alternating knots in $S\times I$.
Pacific J. Math. {\bf 203} (2002), no. 1, 1--22.

\bibitem{Boll} B. Bollob\'{a}s, Modern graph theory. Graduate Texts in Mathematics, {\bf 184}, Springer-Verlag, New York, 1998.

\bibitem{BR1} B. Bollob\'{a}s and O. Riordan,  A polynomial for graphs on orientable surfaces,  Proc. London Math. Soc.  \textbf{83}  (2001),  513-531.

\bibitem{BR} B. Bollob\'{a}s and O. Riordan,  A polynomial of graphs on surfaces,  Math. Ann.  \textbf{323}  (2002),  no. 1, 81-96.


\bibitem{CP} S. Chmutov and I. Pak,  The Kauffman bracket and the Bollobas-Riordan polynomial of ribbon graphs, preprint, {\tt arXiv:math.GT/0404475}.

\bibitem{CP2} S. Chmutov and I. Pak,  The Kauffman bracket of virtual links and the Bollob\'{a}s-Riordan polynomial, {\tt arXiv:math.GT/0609012}, to appear in the Moscow Mathematical Journal.

\bibitem{CRR} P. Cotta-Ramusino and M. Rinaldi, On the algebraic structure of link-diagrams on a $2$-dimensional surface,  Comm. Math. Phys.  {\bf 138}  (1991),  no. 1, 137-173.

\bibitem{Da} O. T. Dasbach, D.~Futer, E.~Kalfagianni, \mbox{X.-S. Lin} and N.~W.~Stoltzfus,  The Jones polynomial and dessins d'enfant, preprint, {\tt arXiv:math.GT/0605571}.

\bibitem{dk} H. A. Dye and L. H.  Kauffman,
Minimal surface representations of virtual knots and links.
Algebr. Geom. Topol. {\bf 5} (2005), 509-535.

\bibitem{EMS} J. Ellis-Monaghan and I. Sarmiento,  A duality relation for the topological Tutte polynomial, talk at the AMS Eastern Section Meeting Special Session on Graph and Matroid Invariants,  Bard College, 10/9/2005. {\tt http://academics.smcvt.edu/jellis-monaghan/\#Research}

\bibitem{Hu}S. Huggett,  On tangles and matroids,  J. Knot Theory Ramifications  {\bf 14}  (2005),  no. 7, 919--929.

\bibitem{hp} J. Hoste and J. H.  Przytycki,
An invariant of dichromatic links.
Proc. Amer. Math. Soc. {\bf 105} (1989), no. 4, 1003-1007.

\bibitem{Ja} F. Jaeger,  Tutte polynomials and link polynomials, Proc. Amer. Math. Soc.  {\bf 103 } (1988),  no. 2, 647-654.

\bibitem{JVW} F. Jaeger, D. Vertigan and D. Welsh, On the computational complexity of the Jones and Tutte polynomials, Math. Proc. Cambridge Philos. Soc., {\bf 108} (1990), 35-53.

\bibitem{jones} V. F. R. Jones,
On knot invariants related to some statistical mechanical models.
Pacific J. Math. {\bf 137} (1989), no. 2, 311-334.

\bibitem{ik} K. Inoue and T. Kaneto,
A Jones type invariant of links in the product space of a surface and the real line.
J. Knot Theory Ramifications {\bf 3} (1994), no. 2, 153-161.

\bibitem{kauf} L. H. Kauffman,  New invariants in the theory of knots.  Amer. Math. Monthly  {\bf 95}  (1988),  no. 3, 195-242.

\bibitem{Lik} W. B. R. Lickorish, An introduction to knot theory. Graduate Texts in Mathematics, {\bf 175}. Springer-Verlag, New York, 1997.

\bibitem{Li} J. Lieberum, Skein modules of links in cylinders over surfaces,  Int. J. Math. Math. Sci.  {\bf 32}  (2002),  no. 9, 515--554.

\bibitem{LM} M. Loebl and I. Moffatt, The chromatic polynomial of fatgraphs and its categorification, preprint, {\tt arXiv:math.CO/0511557}.

\bibitem{m} V. O. Manturov,
Kauffman-like polynomial and curves in 2-surfaces.
J. Knot Theory Ramifications {\bf 12} (2003), no. 8, 1145-1153.

\bibitem{Pr} J. H.  Przytycki,  Skein module of links in a handlebody,  Topology '90 (Columbus, OH, 1990),  315-342, Ohio State Univ. Math. Res. Inst. Publ., 1, de Gruyter, Berlin, 1992.

\bibitem{So} A. D. Sokal, The multivariate Tutte polynomial (alias Potts model) for graphs and matroids, Surveys in Combinatorics, 2005, ed.  Bridget S. Webb (Cambridge University Press, 2005), pp. 173-226.

\bibitem{this} M. B. Thistlethwaite,
A spanning tree expansion of the Jones polynomial.
Topology {\bf 26} (1987), no. 3, 297-309.


\bibitem{Tr} L. Traldi,
A dichromatic polynomial for weighted graphs and link polynomials,
Proc. Amer. Math. Soc. {\bf 106} (1989), no. 1, 279-286.

\bibitem{tu} V. G. Turaev, The Conway and Kauffman modules of a solid torus.  J. Soviet Math.  {\bf 52}  (1990),  no. 1, 2799-2805.

\bibitem{Wo} D. Woodall, Tutte polynomial expansions for 2-separable graphs, Discrete Math., {\bf 247} (2002), 201-213.

\end{thebibliography}
\end{document}